\documentclass[letterpaper, 10 pt, conference]{ieeeconf}
\usepackage{amssymb}
\usepackage{color}
\usepackage{bbm}
\usepackage{amsmath}
\usepackage{graphicx}
\usepackage{mathtools}
\usepackage{multirow}
\usepackage{lineno,hyperref}
\usepackage{url}
\usepackage{cleveref}
\usepackage{subcaption}
\usepackage{bm}
\usepackage{tikz}
\usepackage{algorithm}
\usepackage{algorithmicx}
\usepackage{algpseudocode}
\usepackage{caption}
\usepackage{mathtools}
\mathtoolsset{showonlyrefs,showmanualtags}

\DeclareMathAlphabet{\pazocal}{OMS}{zplm}{m}{n}
\usepackage{graphicx}

\newcommand{\ve}{\boldsymbol}

\usepackage{booktabs, siunitx}

\renewcommand{\geq}{\geqslant}
\renewcommand{\epsilon}{\varepsilon}

\usepackage{pgfplots}
\DeclareUnicodeCharacter{2212}{−}
\usepgfplotslibrary{groupplots,dateplot}
\usetikzlibrary{patterns,shapes.arrows}

\modulolinenumbers[5]

\usepackage{multirow}
\usepackage{color}
\usepackage{paralist}
\usepackage[space]{cite}

\pgfplotsset{compat=1.18}

\IEEEoverridecommandlockouts                        

\title{\LARGE \bf
Gain-Scheduled Data-Enabled Predictive Control: \\ A DeePC Approach for Nonlinear Systems
}

\author{Margarita A. Guerrero, Braghadeesh Lakshminarayanan and Cristian R. Rojas% <-this % stops a space
\thanks{}% <-this % stops a space
\thanks{This work was partially supported by the Swedish Research Council under contract number 2023-05170, and by the Wallenberg AI, Autonomous Systems and Software Program (WASP) funded by the Knut and Alice Wallenberg Foundation. The authors are with the Division of Decision and Control Systems, KTH Royal Institute of Technology, 100 44 Stockholm, Sweden (e-mails: mags3@kth.se, blak@kth.se, crro@kth.se).
}%
}

\begin{document}

\maketitle

%%%%%%%%%%%%%%%%%%%%%%%%%%%%%%%%%%%%%%%%%%%%%%%%%%%%%%%%%%%%%%%%%%%%%%%%%%%%%%%%
% \begin{abstract}
% Estimating optimal inputs for trajectory tracking is central to model predictive control (MPC); the challenge is amplified when the plant is unknown and only input–output data are available. DeePC addresses this in the linear time-invariant setting, %yet many practical plants exhibit pronounced operating-point dependence. 
% yet many real systems exhibit strong operating-point dependence. Building on classical linear parameter-varying control, we introduce DeePC-GS, a novel gain-scheduled extension of data-enabled predictive control for unknown, regime-varying systems. The key idea is to replace the single linear predictor by switching between local Hankel matrices, selected online by a measurable scheduler variable %with a scheduled nonparametric predictor assembled from a bank of local Hankel data sets and selected online via a measurable scheduler
% —thereby uniting classical scheduling theory with identification-free, data-driven MPC. On a nonlinear ship-steering benchmark, DeePC-GS demonstrates the effectiveness of our formulation, outperforming state-of-the-art data-driven MPC while maintaining tractable computation.
% \end{abstract}
\begin{abstract}
Model predictive control is a well established control technology
for trajectory tracking. Its use requires the availability of an
accurate model of the plant, but obtaining such a model is often
time consuming and costly. Data-Enabled Predictive Control (DeePC)
addresses this shortcoming in the linear time-invariant setting, by
skipping the model building step and instead relying directly on
input-output data. Unfortunately, many real systems are nonlinear
and exhibit strong operating-point dependence. Building on classical
linear parameter-varying control, we introduce DeePC-GS, a gain-scheduled
extension of DeePC for unknown, regime-varying systems.
The key idea is to allow DeePC to switch between different local Hankel
matrices—selected online via a measurable scheduling variable—thereby
uniting classical gain scheduling tools with identification-free, data-driven MPC.
We test the effectiveness of our DeePC-GS formulation on a
nonlinear ship-steering benchmark, demonstrating that it outperforms
state-of-the-art data-driven MPC while maintaining tractable computation.
\end{abstract}

% \begin{abstract}
% Estimating optimal inputs for trajectory tracking is central to model predictive control (MPC); the challenge is amplified when the plant is unknown and only input–output data are available. Data-Enabled Predictive Control (DeePC) addresses this in the linear time-invariant setting, yet many real systems exhibit strong operating-point dependence. Building on classical linear parameter-varying control, we introduce DeePC-GS, a gain-scheduled extension of DeePC for unknown, regime-varying systems. The key idea is to replace the single linear predictor by switching between local Hankel matrices—selected online via a measurable scheduling variable—thereby uniting classical scheduling theory with identification-free, data-driven MPC. On a nonlinear ship-steering benchmark, DeePC-GS demonstrates the effectiveness of our formulation, outperforming state-of-the-art data-driven MPC while maintaining tractable computation.
% \end{abstract}


%%%%%%%%%%%%%%%%%%%%%%%%%%%%%%%%%%%%%%%%%%%%%%%%%%%%%%%%%%%%%%%%%%%%%%%%%%%%%%%%
\section{Introduction}

Control engineering has been dominated by model-based design, where traditional system identification procedures for linear systems~\cite{Ljung:99,Soderstrom-Stoica-89} have proven highly effective. These methods, however, often require careful experiment design and can be time consuming. This paves the way for model-free control design procedures.

Several model-free (data-driven) control design procedures have been proposed in the literature~\cite{hjalmarsson1998iterative,campi2002virtual,karimi2004iterative,formentin2012non} that bypass explicit system identification, and controller parameters are tuned directly from input-output measurements from the plant.
%These methods typically assume a fixed controller structure (e.g., PI or PID) and use input-output data from experiments to directly tune the controller parameters to optimize a performance metric, such as tracking error. While effective, this approach is fundamentally about finding the optimal parameters for a pre-defined controller. 

Recently, a novel approach, called the \emph{Data-Enabled Predictive Control (DeePC)} framework~\cite{deepc}, emerged from the behavioral systems theory pioneered by Willems et al.~\cite{willems2005note}. In \cite{willems2005note} it is shown that the behavior of a linear time-invariant (LTI) system can be fully characterized by a sufficiently rich collection of its input-output trajectories. Specifically, \emph{Willems' Fundamental Lemma} states that any new input-output trajectory of an LTI system must lie within the span of previously collected trajectories, provided the input data is \emph{persistently exciting}~\cite{Markovsky_23}. The implication of Willems' lemma is that instead of using data to tune a controller, the data itself can be used as a non-parametric, implicit model to predict and optimize the system's future behavior directly.
%This foundational result has paved the way for model-free control design, most notably
% Persistence of excitation is a condition ensuring the input signal is rich enough to excite all relevant system dynamics, preventing the data from lying in a lower-dimensional subspace.
 
DeePC formulates the predictive control problem directly from input-output data, bypassing the need for explicit model identification. While powerful, the original DeePC formulation is restricted to LTI systems, subsequent works have extended these ideas to nonlinear systems~\cite{huang2023robust,lazar2024basis}, often by employing techniques like basis function expansions or Koopman operators~\cite{koopmanism} to lift the nonlinear dynamics into a higher-dimensional linear space where LTI-based DeePC can be applied. An alternative procedure is to use the existing linear DeePC approaches around several different operating points,  around which the plant behaves almost linearly. Indeed, classical model-based control theory has effectively addressed the challenge of varying operating conditions in nonlinear systems through \emph{gain scheduling}~\cite{astrom1994adaptive,khalil2002nonlinear,Rugh_20}. This well-established technique involves designing a family of linear controllers, each tailored to a specific operating point derived from a linearized model. A \emph{scheduling variable}, which tracks the system's current operating condition in real-time, is then used to interpolate between these controllers. %, effectively adapting the control law as the system's dynamics change.

% However, a common limitation of these data-driven approaches is their reliance on the assumption that the system operates under a single, fixed operating condition. If the underlying nonlinear system experiences significant changes in its operating point, a controller based on a single dataset will inevitably suffer from performance degradation.

% This raises a compelling question: can the adaptive philosophy of gain scheduling be integrated with the model-free benefit of data-driven predictive control?

In this paper, we integrate gain scheduling with Data-Enabled Predictive Control, by proposing a \emph{Gain-Scheduled Data-Enabled Predictive Controller (DeePC-GS)}. Our approach leverages the core principle of gain scheduling not to interpolate controller parameters, but to manage and switch between different datasets, each corresponding to a distinct operating region of the nonlinear system. We select a scheduling variable which triggers a switch to a new data matrix that accurately represents the system's local dynamics. This newly selected dataset is then used within a DeePC formulation to compute the optimal control action. By doing so, we ensure that the predictive controller always relies on data relevant to the current dynamics, overcoming the performance limitations faced by existing data-driven methods for nonlinear systems under varying conditions.

Our main contributions are:
\begin{itemize}
    \item We propose a novel Gain-Scheduled Data-Enabled Predictive Control (DeePC-GS) framework for nonlinear systems.
    \item We present a scheme for switching between local datasets, based on a scheduling variable.
    \item We demonstrate the effectiveness of the proposed DeePC-GS algorithm on a practical example, highlighting its superior performance compared to standard data-driven approaches.
\end{itemize}

The remainder of this paper is organized as follows: In Section~\ref{sec: setup}, we state our problem setup. Section~\ref{sec: prelims} provides some preliminaries on DeePC and gain scheduling. In Section~\ref{sec:approach}, we outline our proposed approach, and demonstrate its efficacy in Section~\ref{sec: Simulation}. Finally, the paper is concluded in Section~\ref{sec: conclusion}.   

\section{Problem Statement} \label{sec: setup}
Consider a nonlinear, discrete-time, time-invariant, multiple-input multiple-output (MIMO) system
\begin{IEEEeqnarray}{rCl}
\label{eq:nlsys}
\begin{aligned}
  \ve{x}(k+1) &=&\ve{f}(\ve{x}(k), \ve{u}(k)) \\
  \ve{y}(k)   &=&\ve{h}(\ve{x}(k), \ve{u}(k)),
\end{aligned}
\end{IEEEeqnarray}
where $\ve{x}(k) \in \mathcal{X} \subset \mathbb{R}^n$ is the internal state, $\ve{u}(k) \in \mathcal{U} \subset \mathbb{R}^m$ is the input into the system, and $\ve{y}(k)\in\mathcal{Y}\subset\mathbb{R}^p$ is the measured output. 
The functions $\ve{f}\colon\mathbb{R}^n\times\mathbb{R}^m\!\to\!\mathbb{R}^n$ and $\ve{h}\colon\mathbb{R}^{n}\times\mathbb{R}^{m}\!\to\!\mathbb{R}^{p}$ govern the state transition and the output mapping, respectively. The objective is to compute an input sequence $\ve{u}_{k:k+N-1}$ that makes $\ve{y}_{k:k+N-1}$ behave in a desired manner---e.g., track a given reference $\ve{r}_{k:k+N-1}$---while satisfying input/output constraints. We stack sequences column-wise by default: for any signal $\ve{s}$, $\ve{s}_{a{:}b}:=\mathrm{col}(\ve{s}_{k+a},\ldots,\ve{s}_{k+b}):=[\,\ve{s}_{k+a}^{\top}\ \cdots\ \ve{s}_{k+b}^{\top}\,]^{\top}$.

A well-known way to achieve this objective, when a model of~\eqref{eq:nlsys} is available, is to solve a receding-horizon nonlinear MPC tracking problem:
% \begin{IEEEeqnarray}{rCl}
% \label{eq:nmpc}
% \begin{aligned}
% \min_{\ve u,\ve x,\ve y}\quad
% & \sum_{j=0}^{N-1}\!\Big(\,\|\ve{y}_j-\ve{r}_{k+j}\|_{Q}^{2}+\|\ve{u}_j\|_{R}^{2}+\|\Delta\ve{u}_j\|_{S}\Big) \\[2pt]
% \text{s.t.}\quad
% & \text{Eq.\,\eqref{eq:nlsys} holds for } j=0,\ldots,N{-}1,\\
% & \ve{x}_0=\hat{\ve{x}}(k),\\
% & \ve{u}_j \in \mathcal{U},\quad j=0,\ldots,N{-}1,\\
% & \ve{y}_j \in \mathcal{Y},\quad j=0,\ldots,N{-}1.
% \end{aligned}
% \end{IEEEeqnarray}
% \begin{IEEEeqnarray}{rCl}
% \label{eq:nmpc}
% \begin{aligned}
% \min_{\ve u,\ve x,\ve y}\quad
% & \sum_{j=0}^{N-1}\!\Big(\,\|\ve{y}_j-\ve{r}_{k+j}\|_{Q}^{2}+\|\ve{u}_j\|_{R}^{2}+\|\Delta\ve{u}_j\|_{S}^{2}\Big) \\[2pt]
% \text{s.t.}\quad
% & \text{Eq.\,\eqref{eq:nlsys} holds for } j=0,\ldots,N{-}1,\\
% & \ve{x}_0=\hat{\ve{x}}(k),\\
% & \big(\,\ve y_{0:N-1}(k),\,\ve u_{0:N-1}(k)\,\big)\in \mathcal{Y}^{N}\times\mathcal{U}^{N}.
% % \quad\text{(optionally add rate limits: } \Delta \ve u_{[0,N-1]}(k)\in \Delta\mathcal U^{N}\text{)}
% \end{aligned}
% \end{IEEEeqnarray}
\begin{IEEEeqnarray}{rCl}
\label{eq:nmpc}
\begin{aligned}
\min_{\ve u,\ve x,\ve y}\quad
& \sum_{j=0}^{N-1}\!\Big(\,\|\ve y_j-\ve r_{k+j}\|_{\ve Q}^{2}+\|\ve u_j\|_{\ve R}^{2}+\|\Delta\ve u_j\|_{\ve S}^{2}\Big) \\[2pt]
\text{s.t.}\quad
& \eqref{eq:nlsys},\quad \ve x_0=\hat{\ve x}(k),\\
& \big(\ve y_{0:N-1}(k),\,\ve u_{0:N-1}(k)\big)\in \mathcal{Y}^{N}\times\mathcal{U}^{N}.
\end{aligned}
\end{IEEEeqnarray}

Here $N\!\in\!\mathbb{N}$ is the prediction horizon, $k$ is the current time index, %$i\!\in\!\{0,\ldots,N{-}1\}$ indexes stages within the horizon, 
$\hat{\ve{x}}(k)$ is the measured/estimated state at time $k$ (based on measurements up to time $k$), $\ve{Q}\succeq \ve{0}$, $\ve{R},\ve{S}\succ \ve{0}$, and $\|\ve v\|_{\ve M}^{2}\!:=\!\ve v^{\!\top} \ve{M} \ve v$ denotes a weighted Euclidean norm.

When $\ve f$ and $\ve h$ are unknown and only input–output data are available, the model-based program~\eqref{eq:nmpc} becomes untenable and a data-driven formulation is required. A widely used approach is \emph{DeePC}, which casts finite-horizon optimal control in MPC form while constructing a nonparametric predictor from Hankel matrices of persistently exciting data—thus avoiding explicit system identification. However, for strongly nonlinear plants with multiple equilibria or operating regions, a single offline Hankel matrix built from mixed trajectories may fail to be locally representative; prediction quality degrades because the LTI assumptions underlying the DeePC predictor generally do not hold. Motivated by these limitations, we introduce in Section~\ref{sec:approach} \emph{DeePC-GS}, a gain-scheduled, fully data-driven predictive control scheme that leverages LTI DeePC predictors in a nonlinear setting by switching between local Hankel models via a measurable scheduler, thereby retaining linear MPC machinery while switching between different operating regions.

\section{Preliminaries} \label{sec: prelims}
This section reviews the theoretical underpinnings of DeePC and gain scheduling on which DeePC-GS is built. The treatment is self-contained and draws on \cite{deepc,Rugh_20}.

\subsection{DeePC}
Consider a discrete-time LTI system with realization $(\ve{A},\ve{B},\ve{C},\ve{D})$, input $\ve u\in\mathbb{R}^m$, output $\ve y\in\mathbb{R}^p$, and state $\ve x\in\mathbb{R}^n$. The behavior restricted to a window of length-$L$ trajectories is the set of pairs $(\ve u,\ve y)$ for which there exists a state sequence satisfying the state/output equations, namely
$\mathcal{B}(\ve{A},\ve{B},\ve{C},\ve{D})=\big\{\mathrm{col}(\ve u,\ve y)\ \big|\ \exists\,\ve x\ \text{s.t.}\ (q\ve x)_k=\ve A\ve x+\ve B\ve u,\ \ve y=\ve C\ve x+\ve D\ve u\big\}$---$q$ denotes the forward-shift operator.
For an observable pair $(\ve A,\ve C)$, observability over windows is characterized by the extended observability matrix $\mathcal O_\ell(\ve A,\ve C)=\mathrm{col}(\ve C,\ve {CA},\dots,\ve{CA}^{\ell-1})$, with the \emph{lag} $\ell$ defined as the smallest integer for which $\mathrm{rank}\,\mathcal O_\ell=n$.

From an \emph{offline} noiseless experiment on the plant, let $\ve u^{d}=\mathrm{col}(\ve u^{d}_{1},\ldots,\ve u^{d}_{T})\in\mathbb{R}^{mT}$ and $\ve y^{d}=\mathrm{col}(\ve y^{d}_{1},\ldots,\ve y^{d}_{T})\in\mathbb{R}^{pT}$ denote the collected sequences (superscript $d$ marks offline data). Fix the window depth $L:=T_{\mathrm{ini}}+N$, where $T_{\mathrm{ini}}$ counts the most recent samples used to match online measurements (thus implicitly initializing the unknown state) and $N$ is the future horizon over which we predict/optimize. Define the depth-$L$ Hankel matrices $\ve{\mathcal{H}}_L(\ve u^{d})\in\mathbb{R}^{mL\times H_c}$ and $\ve{\mathcal{H}}_L(\ve y^{d})\in\mathbb{R}^{pL\times H_c}$—with $H_c:=T{-}L{+}1$—whose columns are time-shifted length-$L$ windows of $\ve u^{d}$ and $\ve y^{d}$. If $\ve u^d$ is \emph{persistently exciting (PE)} of order $L{+}n$—i.e., the depth-$(L{+}n)$ Hankel $\ve{\mathcal{H}}_{L{+}n}(\ve u^d)$ has full row rank $m(L{+}n)$ (a sufficient condition is $T\!\ge\!(m{+}1)(L{+}n){-}1$)—then for any $T_{\mathrm{ini}}\!\geq\!\ell$ the column space of the stacked matrix $[\,\ve{\mathcal{H}}_L(\ve u^{d});\,\ve{\mathcal{H}}_L(\ve y^{d})\,]$ coincides with the set of all length-$L$ trajectories; equivalently, its rank is $mL{+}n$~\cite{Markovsky_23,deepc}. %\emph{Note that the excitation condition is stated at order $L{+}n$, whereas the DeePC predictor itself uses depth-$L$ Hankel matrices; the ``$+n$'' pertains to the richness requirement that absorbs the unknown initial state, not to the matrices used in the optimization.}

We exploit this fact to build a purely data-driven predictor by partitioning the depth-$L$ Hankel matrices into past/future blocks via $\ve U_P:=\ve{\mathcal{H}}_{T_{\mathrm{ini}}}(\ve u^{d})$, $\ve Y_P:=\ve{\mathcal{H}}_{T_{\mathrm{ini}}}(\ve y^{d})$, $\ve U_F:=\ve{\mathcal{H}}_{N}(\ve u^{d})$, $\ve Y_F:=\ve{\mathcal{H}}_{N}(\ve y^{d})$ %, where $H_c:=T{-}L{+}1$ is the number of columns (i.e., available length-$L$ windows) 
with $\ve U_P\!\in\!\mathbb{R}^{mT_{\mathrm{ini}}\times H_c}$, $\ve Y_P\!\in\!\mathbb{R}^{pT_{\mathrm{ini}}\times H_c}$, $\ve U_F\!\in\!\mathbb{R}^{mN\times H_c}$, $\ve Y_F\!\in\!\mathbb{R}^{pN\times H_c}$. Given online past windows $\ve u_{\mathrm{ini}}=\ve u_{\,k-T_{\mathrm{ini}}:k-1}$ and $\ve  y_{\mathrm{ini}}=\ve y_{\,k-T_{\mathrm{ini}}+1:k}$, and future sequences $(\ve u,\ve y)$ of length $N$, DeePC enforces the data equation
% \begin{IEEEeqnarray}{rCl}
% \label{eq:deepc_data}
% \centering
% \begin{bmatrix}
% U_P\\ Y_P\\ U_F\\ Y_F
% \end{bmatrix}\ve g
% =
% \begin{bmatrix}
% \ve u_{\mathrm{ini}}\\ \ve y_{\mathrm{ini}}\\ \ve u\\ \ve y
% \end{bmatrix},\qquad \ve g\in\mathbb{R}^{H_c}.
% \end{IEEEeqnarray}
\begin{IEEEeqnarray}{rCl}
\label{eq:deepc_data}
\operatorname{col}(\ve U_P,\ve Y_P,\ve U_F,\ve Y_F)\,\ve g
&=&
\operatorname{col}(\ve u_{\mathrm{ini}},\ve y_{\mathrm{ini}},\ve u,\ve y),
\end{IEEEeqnarray}
The past block $(\ve U_P,\ve Y_P)$ enforces behavioral consistency with $(\ve u_{\mathrm{ini}},\ve y_{\mathrm{ini}})$, implicitly selecting the internal state; the future block $(\ve U_F,\ve Y_F)$ maps $\ve g \in \mathbb{R}^{H_c}$ to the predicted $(\ve u,\ve y)$. We refer to the \emph{DeePC predictor} as the full data equation built from $(\ve U_P,\ve Y_P,\ve U_F,\ve Y_F)$; cf.~\eqref{eq:deepc_data}.

DeePC solves, in a receding horizon fashion, an MPC program of the same structure as our \eqref{eq:nlsys}-based design—minimizing a tracking/effort stage cost (e.g., $\|\ve y-\ve r\|_{\ve Q}^{2}+\ell(\ve u)$) under input/output constraints—but with the model-dynamics constraint replaced by the data equation \eqref{eq:deepc_data}. With noisy data, this constraint is imposed softly via slack variables and regularization on $\ve g$ and on the past residuals. In Section~\ref{sec:gs} we briefly recall classical gain scheduling and integrate it with DeePC, yielding DeePC-GS.

\subsection{Gain Scheduling}\label{sec:gs}
Real plants rarely behave linearly across their operating envelopes: aerodynamic/thermal effects, actuator limits, loads, and set-point–dependent physics all make the dynamics vary with measurable operating conditions. Gain scheduling addresses this by precomputing a \emph{family of linearizations} and a matching \emph{family of linear controllers}, both indexed by a measurable scheduling variable $\rho(k)$. The controller parameters then vary online with $\rho(k)$, retaining linear design intuition while covering a range of operating conditions. Mathematically, the gain-scheduling methodology adapts a controller to the changing dynamics of the nonlinear system \eqref{eq:nlsys} by following a three-step design process.

First, the nonlinear system is linearized around a family of equilibrium (steady-state) points $(\ve{x}_{\mathrm{ss}}(\rho),\,\ve{u}_{\mathrm{ss}}(\rho))$ parameterized by the scheduling variable $\rho$, yielding an Linear Parameter-Varying (LPV) model in deviation coordinates: \begin{align}
\ve{x}_\delta(k{+}1)&=\ve{A}(\rho)\,\ve{x}_\delta(k)+\ve{B}(\rho)\,\ve{u}_\delta(k),\\
\ve{y}_\delta(k)&=\ve{C}(\rho)\,\ve{x}_\delta(k)+\ve{D}(\rho)\,\ve{u}_\delta(k).
\end{align}
Here $\ve{A}(\rho),\ve{B}(\rho),\ve{C}(\rho),\ve{D}(\rho)$ are the Jacobians of $\ve f$ and $\ve h$ evaluated at the equilibrium point associated with $\rho$, with $\rho$ treated as frozen (constant)—hence a frozen-parameter LTI plant at each operating condition.

Second, for this LPV surrogate, a family of linear controllers is designed—one for each frozen value of $\rho$. The controller parameters (e.g., gain matrices $(\ve{F}(\rho),\ve{G}(\rho),\ve{L}(\rho),\ve{M}(\rho))$) are thus functions of $\rho$; designs may be computed continuously in $\rho$ or at a finite grid of operating points with interpolation, ensuring consistency with the equilibrium map $(\ve{x}_{\mathrm{ss}}(\rho),\ve{u}_{\mathrm{ss}}(\rho))$ and recovery of the intended linear closed loop at each frozen $\rho$.

Finally, a single nonlinear controller is implemented by replacing the static parameter with the real-time scheduling signal $\rho(k)$, producing the scheduled law in state-space form
% \begin{align}
% \ve{z}(k{+}1)&=F(\rho(k))\,\ve{z}(k)+G(\rho(k))\big(\ve{y}(k)-\ve{y}_{\mathrm{ss}}(\rho(k))\big),\\
% \ve{u}(k)&=L(\rho(k))\,\ve{z}(k)+M(\rho(k))\big(\ve{y}(k)-\ve{y}_{\mathrm{ss}}(\rho(k))\big)\notag\\
% &\quad+\ve{u}_{\mathrm{ss}}(\rho(k)).
% \end{align}
\begin{align}
\ve z(k{+}1)&=\ve{F}(\rho(k))\,\ve z(k)+\ve{G}(\rho(k))\,\ve e_y(k),\\
\ve u(k)&=\ve{L}(\rho(k))\,\ve z(k)+\ve M(\rho(k))\,\ve e_y(k)+\ve u_{\mathrm{ss}}(\rho(k)),
\end{align}
where $\ve e_y(k):=\ve y(k)-\ve y_{\mathrm{ss}}(\rho(k))$, and $\{\ve{F},\ve G,\ve L,\ve M,\ve u_{\mathrm{ss}},\ve y_{\mathrm{ss}}\}$ are scheduled by $\rho(k)$.
% where $\ve e_y(k):=\ve y(k)-\ve y_{\mathrm{ss}}(\rho(k))$. The map $\rho\mapsto\{F,G,L,M,\ve{u}_{\mathrm{ss}},\ve{y}_{\mathrm{ss}}\}$ then adapts the control action as operating conditions change. 

Our scheduled variant, DeePC-GS, retains this template and applies scheduling to the predictor $[\,\ve{\mathcal{H}}_L(\ve u^{d});\,\ve{\mathcal{H}}_L(\ve y^{d})\,].$

\section{Proposed Approach}\label{sec:approach}

In this section, we propose \emph{DeePC-GS}, a gain-scheduled, fully data-driven predictive control scheme. While standard DeePC implicitly relies on a single LTI predictor, DeePC-GS constructs a \emph{scheduled nonparametric predictor} from a bank of local input--output datasets indexed by a measurable scheduling variable $\rho(k)$. %We first state the certainty-equivalence MPC problem using this scheduled predictor and then present the \emph{region-selection} realization.\\

\subsection{Data Collection}
We consider an unknown, discrete-time, nonlinear plant $ \mathcal{P} $ with $ m $ inputs and $ p $ outputs. Let $ \mathcal{B}_T(\mathcal{P}) $ denote the set of all length-$ T $ admissible input--output trajectories of $ \mathcal{P} $. The plant operates near a finite family of operating regions, indexed by $ i\in \mathcal{I}=\{1,\dots,M\}\subset\mathbb{N}$ and parameterized by one or more measurable scheduling variables $ \rho(k) $. We write $ \mathcal{B}_T^{(i)}(\mathcal{P}) $ for the restriction of $ \mathcal{B}_T(\mathcal{P}) $ to region $ i $, and assume that the local I/O behavior is well approximated by a minimal LTI surrogate $ \ve{\Sigma}_i=(\ve{A}_i,\ve{B}_i,\ve{C}_i,\ve{D}_i) $ of order $ n_i $, so that $ \mathcal{B}_T^{(i)}(\mathcal{P}) \approx \mathcal{B}_T(\ve \Sigma_i) $; hence Willems’ Fundamental Lemma applies~\cite{willems2005note}.%---this assumption does not turn $\mathcal{P}$ into an LTI system; it only states that locally its behavior can be captured by $\ve{\Sigma}_i$, which is the object for which Willems' Lemma applies.

Offline, we collect $M$ local input--output datasets around distinct operating points parameterized by the scheduling variable $\rho$. Let $\pi\colon \mathcal{D}_\rho\to\mathcal{I}$ be a (measurable) selection map from the scheduler domain $\mathcal{D}_\rho$ to the index set $\mathcal{I}$, and define $i_k := \pi(\rho(k))$. For each $i\in\mathcal{I}$, let $(\ve{u}_i^{d},\ve{y}_i^{d}) \in\mathbb{R}^{mT_i}\times\mathbb{R}^{pT_i}$ denote \emph{offline} I/O sequences of total length $T_i$---i.e., $\ve u_i^{d}=\mathrm{col}\big(\ve u_i^{d}(0),\ldots,\ve u_i^{d}(T_i{-}1)\big),
\ve y_i^{d}=\mathrm{col}\big(\ve y_i^{d}(0),\ldots,\ve y_i^{d}(T_i{-}1)\big)$---measured from the unknown plant $\mathcal{P}$ under inputs that are persistently exciting of order $L+n_i$, where $L:=T_{\mathrm{ini}}+N$ and $n_i := n\big(\mathcal{B}_T^{(i)}(\mathcal{P})\big)$ denotes the local (behavioral) order. A sufficient dataset length is $T_i \ge (m+1)(L+n_i)-1$. The Hankel columns are built from contiguous windows of length $L$ on which the region label is constant---i.e., $i_k=i$ for all $k$ in the window. From these sequences we form the Hankel partitions
% \[
% \begin{bmatrix}
% U_{P,i}\\[1pt] Y_{P,i}\\[1pt] U_{F,i}\\[1pt] Y_{F,i}
% \end{bmatrix}
% :=
% \begin{bmatrix}
% H_{T_{\mathrm{ini}}}(\ve{u}_i^{d})\\[1pt]
% H_{T_{\mathrm{ini}}}(\ve{y}_i^{d})\\[1pt]
% H_{N}(\ve{u}_i^{d})\\[1pt]
% H_{N}(\ve{y}_i^{d})
% \end{bmatrix}.
% \]
\[
\begin{aligned}
\mathrm{col}(\ve U_{P,i},& \ve Y_{P,i},\ve U_{F,i},\ve Y_{F,i}) := \\
&\mathrm{col}\big(\ve{\mathcal{H}}_{T_{\mathrm{ini}}}(\ve u_i^{d}),\,\ve{\mathcal{H}}_{T_{\mathrm{ini}}}(\ve y_i^{d}),\,\ve{\mathcal{H}}_{N}(\ve u_i^{d}),\,\ve{\mathcal{H}}_{N}(\ve y_i^{d})\big).
\end{aligned}
\]
% \[
% \begin{aligned}
% [\,U_{P,i};\ & Y_{P,i};\ U_{F,i};\ Y_{F,i}\,]\; \coloneqq\;  \\[-2pt]
% &[\,H_{T_{\mathrm{ini}}}(\ve u_i^{d});\ H_{T_{\mathrm{ini}}}(\ve y_i^{d});\ H_{N}(\ve u_i^{d});\ H_{N}(\ve y_i^{d})\,].
% \end{aligned}
% \]

Intuitively, each dataset $i$ yields a local nonparametric predictor valid near its operating point.

\subsection{DeePC-GS Optimization}\label{sec:Optimization}
Given past windows $(\ve{u}_{\mathrm{ini}}, \ve{y}_{\mathrm{ini}})$ and the selected region $i_k$, we solve at each step $k$ the following receding-horizon program:
% \begin{IEEEeqnarray}{rCl}
% \label{eq:deepc_gs_region}
% \begin{aligned}
% \min_{\ve{u},\,\ve{y},\,\ve{g},\,\ve{s}_y,\,\ve{s}_u}\;\; &
% \sum_{j=0}^{N-1}\!\Big(\,\|\ve{y}_j-\ve{r}_{k+j}\|^2_{Q_{i_k}}
% \;+\; \ell_{i_k}(\ve{u}_j)\Big)\;
% +
% \\[2pt]
% &\quad
% \lambda_{g,i_k}\|\ve{g}\|_{1}
% +\lambda_{y,i_k}\|\ve{s}_y\|_{1}
% +\lambda_{u,i_k}\|\ve{s}_u\|_{1}
% \\[4pt]
% \text{s.t. }\; &
% \begin{bmatrix}
% U_{P,{i_k}}\\[2pt] Y_{P,{i_k}}\\[2pt] U_{F,{i_k}}\\[2pt] Y_{F,{i_k}}
% \end{bmatrix}\ve{g}
% =
% \begin{bmatrix}
% \ve{u}_{\mathrm{ini}}+\ve{s}_u\\[2pt] \ve{y}_{\mathrm{ini}}+\ve{s}_y\\[2pt] \ve{u}\\[2pt] \ve{y}
% \end{bmatrix},
% \\[4pt]
% & \ve{u}_j \in \mathcal{U}\,,\quad j=0,\ldots,N-1,\\[2pt]
% & \ve{y}_j \in \mathcal{Y}\,,\quad j=0,\ldots,N-1\,.
% \end{aligned}
% \end{IEEEeqnarray}
\begin{IEEEeqnarray}{rCl}
\label{eq:deepc_gs_region}
\begin{aligned}
\min_{\ve{u},\,\ve{y},\,\ve{g},\,\ve{s}_y,\,\ve{s}_u}&\;\; 
\sum_{j=0}^{N-1}\!\Big(\,\|\ve{y}_j-\ve{r}_{k+j}\|^2_{\ve{Q}_{i_k}}+\ell_{i_k}(\ve{u}_j)\Big)\\
&\;\;+\;\lambda_{g,i_k}\|\ve{g}\|_{1}
+\lambda_{y,i_k}\|\ve{s}_y\|_{1}
+\lambda_{u,i_k}\|\ve{s}_u\|_{1} \\[7pt]
\text{s.t. }\quad & \; \;
\operatorname{col}(\ve U_{P,i_k},\ve Y_{P,i_k},\ve U_{F,i_k},\ve Y_{F,i_k})\,\ve g \\[-1pt]
&\!\hphantom{=}=\;\operatorname{col}(\ve u_{\mathrm{ini}}+\ve s_u,\;\ve y_{\mathrm{ini}}+\ve s_y,\;\ve u,\;\ve y),\\[2pt]
& \;\;\;\big(\ve y_{0:N-1}(k),\,\ve u_{0:N-1}(k)\big)\in \mathcal{Y}^{N}\times\mathcal{U}^{N}.
\end{aligned}
\end{IEEEeqnarray}
Both the tracking term and the stage cost $\ell\colon \mathbb{R}^{mN} \to \mathbb{R}$ are region-dependent via $(\ve{Q}_{i_k},\ve{R}_{i_k},\ve{S}_{i_k})$, enabling operating-point–specific tracking/effort/rate trade-offs; in particular, $\ell_{i_k}(\cdot)$ penalizes input magnitude and rate via $\ve{R}_{i_k}$ and $\ve{S}_{i_k}$. An $\ell_1$ penalty on $\ve g$ stabilizes the optimization and promotes sparsity (selecting few Hankel columns), while measurement noise and local surrogate mismatch—due to nonlinearities and $\rho$ variations over the horizon, including intra-/inter-region transitions—are absorbed by the slack variables $\ve s_y,\ve s_u$ weighted by $(\lambda_{y,i_k},\lambda_{u,i_k})$.
%The tracking error term and the stage cost $\ell\colon \mathbb{R}^{mN} \to \mathbb{R}$ are region dependent through  $(Q_{i_k},R_{i_k},S_{i_k})$, allowing different tracking–effort–rate trade-offs per operating region; the latter penalizes control actions such as input magnitude $\|\ve{u}_k\|_{R_{i_k}}$ and control-rate $\|\Delta\ve{u}_k\|_{S_{i_k}}$. %The positive semi-definite matrix $\ve{Q} \in \mathbb{S}^{pN\times pN}_{\geq0}$ and the positive definite matrices $\ve{R}, \textbf{ }\ve{S} \in  \mathbb{S}^{mN\times mN}_{>0}$ are the cost matrices. 
%Regularization on $\ve{g}$ primarily improves numerical conditioning while robustness to measurement noise and local surrogate mismatch---arising from plant nonlinearities and $ \rho $ variations over the prediction horizon inducing intra-/inter-region transitions-- is handled via the soft relaxations $\ve{s}_y,\ve{s}_u$, whose magnitudes are penalized by $\lambda_{y,i_k}$ and $\lambda_{u,i_k}$.
\vspace{1.2 em}
\subsection{Computational footprint}
At each step we solve \eqref{eq:deepc_gs_region} for the active region $i_k$. Region~$i$ uses an offline slice of length $T_i$, yielding $n_{g,i}=T_i{-}L{+}1$ Hankel columns (the dimension of the coefficient vector $g$). With a dense interior--point method, the per--iteration work is dominated by a dense Cholesky/LDL$^\top$ factorization of a KKT (or normal--equation) matrix of size $n_{g,i_k}\!\times n_{g,i_k}$, which amounts to $O(n_{g,i_k}^3)$ flops per iteration (in our configurations $n_{g,i}\gg (m{+}p)N$, so the $(m{+}p)N$ block is negligible). Under a fixed total data budget $T_{\mathrm{tot}}=\sum_i T_i$, increasing the number of regions $M$ typically shortens each slice---e.g., $T_i\approx T_{\mathrm{tot}}/M$ if split evenly---so $n_{g,i}$ shrinks and the dominant cost scales like $O(((T_{\mathrm{tot}}/M)-L+1)^3)$. %In practice, $T_i$ must still satisfy the PE lower bound $T_i\ge (m+1)(L+n_i)-1$, which limits how large $M$ can be.

For slow processes and long horizons, the Hankel partitions can become large and ill--conditioned. Without changing \eqref{eq:deepc_gs_region}, we can reduce complexity by splitting each region's data bank into a few slices (e.g., by the regulated output's trend---rising/steady/falling) and, online, instantiating the predictor using only the Hankel columns from the matching slice. This reduces the effective value of $n_{g,i}$, improves conditioning, and further lowers computational cost while preserving locality, provided the slice still satisfies the PE length requirement; see Sec.~\ref{sec:expset} for details.

%For slow processes and long horizons, the Hankel partitions can become large and ill-conditioned. Without changing the optimization in \eqref{eq:deepc_gs_region}, we reduce complexity by partitioning—within each region $i$—the local data bank into a few specific subsets ("slices") defined by an operational condition (e.g., the trend of the regulated output—rising/holding/falling). At run time, we instantiate the predictor for region $i$ using only the Hankel columns from the slice consistent with the current condition. This column restriction reduces the effective number of columns, improves numerical conditioning and computational cost, and preserves locality. Full details of this heuristic are provided in Section~\ref{sec:expset}.
%; these slacks preserve feasibility when measurements are noisy or the local surrogate is imperfect.

\subsection{Region switching (hard vs.\ soft)}\label{sec:region-switching}
At runtime, region updates are driven by $i_k=\pi(\rho(k))$.

\noindent\textbf{Hard switching.} A dataset switch is triggered by a hysteresis condition $\mathrm{HYST}_\rho$ on $\rho$—a symmetric deadband around each boundary: we commit a change only if $\rho$ exceeds the nearest boundary by $\mathrm{HYST}_\rho$, preventing chatter. Upon trigger, the DeePC predictor is immediately re-indexed—i.e., we switch to the Hankel partitions $\ve U_{P,i_k},\ve Y_{P,i_k},\ve U_{F,i_k},\ve Y_{F,i_k}$—adopting region-specific bounds when applicable. Since the new dataset may induce a sudden control move, we apply a short \emph{temporal input blend} of duration $T_{\mathrm{BLEND}}$: if the change $i_{\mathrm{old}}\!\to i_{\mathrm{new}}$ is detected at sample $k_0$, then for $k\ge k_0$ and sampling time $\Delta t>0$, the applied input is
\vspace{-0.8mm}
\begin{IEEEeqnarray}{rCl}
\label{eq:hardble}
\begin{aligned}
\ve u_{\mathrm{applied}}(k) &= (1-\alpha_k)\,\ve u^{\mathrm{old}}(k)+\alpha_k\,\ve u^{\mathrm{new}}(k),
% \alpha_k &= \phi\!\Big(\frac{(k - k_0)\Delta t}{T_{\mathrm{BLEND}}}\Big),
\end{aligned}\vspace{-0.5mm}
\end{IEEEeqnarray}
where $\alpha_k=\phi\!\left(\tfrac{(k-k_0)\Delta t}{T_{\mathrm{BLEND}}}\right)$  and $\phi$ is the $C^2$ function $\phi(s)=0$ for $s\le0$, $= 10s^{3}-15s^{4}+6s^{5}$ for $0<s<1$, and $=1$ for $s\ge1$, yielding a jerk-limited transition.
%with $s_k=\frac{(k - k_0)\Delta t}{T_{\mathrm{BLEND}}}$ and the \emph{smootherstep} interpolation $\phi(s)$, defined by $\phi(s)=0$ for $s\le 0$, $\phi(s)=10s^{3}-15s^{4}+6s^{5}$ for $0<s<1$, and $\phi(s)=1$ for $s\ge 1$,
%which is $C^2$ and has zero slope at the endpoints, yielding a jerk-limited transition. %During blending we enforce feasibility by applying the optimization over a conservative input set (e.g., $\mathcal{U}_{i_{\mathrm{old}}}\cap\mathcal{U}_{i_{\mathrm{new}}}$) so that the convex combination remains admissible. 
%We also include \emph{reverse-dwell} abort/commit logic if $\rho$ returns across the switching boundary for at least $T_{REV}$.
We include \emph{reverse-dwell} abort logic when $\rho$ returns to its previous region.

\noindent\textbf{Soft switching.} We pre-smooth region-dependent hyperparameters as a function of $\rho(k)$ in a symmetric overlap band of half-width $\mathrm{\rho_{BAND}}$ around each boundary.  When $\rho$ lies in that band toward the upper neighbor (resp.\ lower), we form the normalized progress $s$ and the weight $\beta=\phi(s)$ via
\vspace{-1mm} % o -1mm, -3mm según necesites
\begin{IEEEeqnarray}{rCl}\label{eq:smooth}
s=\big(\rho(k)-(\mathrm{hi/low}-\rho_{\mathrm{BAND}})\big)/(2\rho_{\mathrm{BAND}}).
\vspace{-1mm} % o -1mm, -3mm según necesites
\end{IEEEeqnarray}
% \begin{IEEEeqnarray}{rCl}\label{eq:smooth}
% s=\frac{\rho(k)-(\mathrm{hi/low}-\mathrm{\rho_{BAND}})}{2\,\mathrm{\rho_{BAND}}},
% \end{IEEEeqnarray}
and blend the weights by a convex combination, $
K^{\mathrm{eff}}=(1-\beta)\,K_{i_{\mathrm{old}}}+\beta\,K_{i_{\mathrm{new}}}$, where $K\in\big\{\ve{Q},\,\ve{R},\,\ve{S},\,\,\lambda_{g},\,\lambda_{y},\,\lambda_{u}\big\}$. This $\rho$-based soft switching reduces abrupt changes of the optimization landscape near region boundaries and complements the temporal input blend.

\subsection{DeePC-GS Algorithm}
%We now present the DeePC-GS \emph{gain-scheduled, data-driven predictive controller} in Algorithm~\ref{alg:deepcgs}. At each sampling instant, we apply only the first control move of the $N$-step plan computed from \eqref{eq:deepc_gs_region}, then shift the past windows and re-solve at $k{+}1$ with updated measurements—i.e., a standard receding-horizon MPC implementation.
We now present the \emph{gain-scheduled, Data-Enabled Predictive Controller} DeePC-GS in Algorithm~\ref{alg:deepcgs}. At each sampling instant, we solve \eqref{eq:deepc_gs_region}, apply only the first control move of the $N$-step plan, and update the DeePC windows by shifting one step—namely, for $k+1$ set $\ve{u}_{\mathrm{ini}}^{+}=\operatorname{col}(\ve u_{k-T_{\mathrm{ini}}+1},\dots,\ve u_{k})$ and $\ve{y}_{\mathrm{ini}}^{+}=\operatorname{col}(\ve y_{k-T_{\mathrm{ini}}+2},\dots,\ve y_{k+1})$ with $\ve u_k$ the applied input and $\ve y_{k+1}$ the newly measured output—then re-solve at $k{+}1$, i.e., a standard receding-horizon MPC implementation.

In Algorithm~\ref{alg:deepcgs}, \texttt{NEIGHBORINBAND} returns the normalized progress $s$ in \eqref{eq:smooth} used to compute $\beta$ when $\rho(k)$ lies within an overlap band at a region boundary; \texttt{SOFTBLENDHYPERPARAMS} soft-switches the region-dependent hyperparameters (Sec.~IV-C); and \texttt{SOLVEDEEPC} solves the DeePC-GS program \eqref{eq:deepc_gs_region}.

\begin{algorithm}[h!]
\footnotesize
\caption{DeePC-GS Algorithm}
\label{alg:deepcgs}
\begin{algorithmic}[1]
\Require Data bank $\{\ve U_{P,i},\ve Y_{P,i},\ve U_{F,i},\ve Y_{F,i}\}_{i\in\mathcal{I}}$; selection map $\pi\colon\mathcal{D}_\rho\!\to\!\mathcal{I}$ and initial $\rho(0)$; initial windows $(\ve u_{\mathrm{ini}},\ve y_{\mathrm{ini}})$; overlap band $\mathrm{\rho_{BAND}}$; hysteresis $\mathrm{HYST}_\rho$; blend time $T_{\mathrm{BLEND}}$.
\State $i_{\text{cur}} \gets \pi(\rho(0))$;$\quad \ve{\mathcal{H}}_{\text{cur}} \gets (\ve U_{P,i_{\text{cur}}},\ve Y_{P,i_{\text{cur}}},\ve U_{F,i_{\text{cur}}},\ve Y_{F,i_{\text{cur}}})$
\State $\texttt{switch} \gets \textbf{Inactive}$
\For{$k=0,1,2,\dots$}
    \State $(\text{overlap},\,i_{\text{new}},\,s)\gets\textsc{NeighborInBand}(i_{\text{cur}},\rho(k),\mathrm{\rho_{BAND}})$
    \If{$\text{overlap}$} \State \textsc{SoftBlendHyperparams}$(i_{\text{cur}},i_{\text{new}},s)$ \Comment{$s$ as in \eqref{eq:smooth}} \EndIf
    \State $\ve u^{\text{cur}} \gets \textsc{SolveDeePC}\big(\ve{\mathcal{H}}_{\text{cur}},(\ve u_{\mathrm{ini}},\ve y_{\mathrm{ini}})\big)$
    \If{$\rho(k)$ \textbf{crosses a boundary (with} $\mathrm{HYST}_\rho$\textbf{)}}
        \State $\texttt{switch} \gets \textbf{Active}$
        \State $\ve{\mathcal{H}}_{\text{new}} \gets (\ve U_{P,i_{\text{new}}},\ve Y_{P,i_{\text{new}}},\ve U_{F,i_{\text{new}}},\ve Y_{F,i_{\text{new}}})$
        \State $\ve u^{\text{new}} \gets \textsc{SolveDeePC}\big(\ve{\mathcal{H}}_{\text{new}},(\ve u_{\mathrm{ini}},\ve y_{\mathrm{ini}})\big)$
        \State $ \ve u_k \gets (1-\alpha_k)\,\ve u^{\text{cur}}[0] + \alpha_k\,\ve u^{\text{new}}[0]$ \Comment{$\alpha_k$ as in \ref{eq:hardble}}
    \Else
        \State $\ve u_k \gets \ve u^{\text{cur}}[0]$
    \EndIf
    \If{$\texttt{switch}=\textbf{Active}$ \textbf{for} $T_{\mathrm{BLEND}}$}
        \State $\texttt{switch}\gets\textbf{Inactive}$;\; $i_{\text{cur}}\gets i_{\text{new}}$;\; $\ve{\mathcal{H}}_{\text{cur}}\gets \ve{\mathcal{H}}_{\text{new}}$
    \EndIf
    \State Apply $\ve u_k$; acquire $\ve y_{k+1}$
    \State $\ve u_{\mathrm{ini}}\gets \operatorname{col}(\ve u_{k-T_{\mathrm{ini}}+1},\ldots,\ve u_k)$
    \State $\ve y_{\mathrm{ini}}\gets \operatorname{col}(\ve y_{k-T_{\mathrm{ini}}+2},\ldots,\ve y_{k+1})$
\EndFor
\end{algorithmic}
\end{algorithm}
\vspace{-1.5em}

\section{Experiments} \label{sec: Simulation}
This section presents simulation studies assessing the performance of DeePC-GS (Algorithm~\ref{alg:deepcgs}) against RoKDeePC~\cite{RoKDeepc}, Koopman-MPC~\cite{koopmanism}, and regularized DeePC~\cite{deepc} on a nonlinear cargo ship--steering benchmark, targeting heading regulation under input constraints.

\subsection{Case Study: Ship Steering Autopilot}
Closed-loop heading control of surface vessels is strongly regime-dependent due to yaw–sway coupling, speed-dependent hydrodynamic derivatives, and environmental loads, causing controllers tuned at a single operating point to degrade off-design. Large cargo vessels exemplify this: long time constants, strong coupling, and a marked dependence on surge speed $U$. These features make DeePC-GS a natural solution: the measurable, slowly varying surge speed $U$—with strong dynamic dependence—serves as an effective scheduling variable to switch \emph{online} between local Hankel predictors \cite{Rugh_20,khalil2002nonlinear}. We therefore adopt a MARINER--hull cargo ship as our nonlinear benchmark \cite{chislett1965planar,fossen1994guidance}.

We model the ship’s horizontal motion with a Nomoto-type maneuvering model augmented with surge (propulsion) dynamics. In a body-fixed frame, the state is $\ve x=[\,v,\ r,\ \psi,\ U\,]^\top$, with sway velocity $v$, yaw rate $r$, heading $\psi$, and surge speed $U$; the control input is $\ve u=[\,\delta,\ \tau\,]^\top$, where $\delta$ is the rudder angle and $\tau$ the propeller thrust. The yaw–sway and surge dynamics are given by
\begin{IEEEeqnarray}{rCl}
\label{eq:surge}
\begin{aligned}
\dot v &= \frac{U}{l}\,a_{11}(U)\,v \;+\; U\,a_{12}(U)\,r \;+\; \frac{U^{2}}{l}\,b_{1}(U)\,\delta, \\
\dot r &= \frac{U}{l^{2}}\,a_{21}(U)\,v \;+\; \frac{U^{2}}{l}\,a_{22}(U)\,r \;+\; \frac{U^{2}}{l^{2}}\,b_{2}(U)\,\delta,\\
\dot \psi &= r, \qquad \dot U = \frac{-X_u\,U + \tau}{m}. %\\
%\dot U &= \frac{-X_u\,U + \tau}{m}. 
\end{aligned}
\end{IEEEeqnarray}
Here $l$ is a characteristic length ($161$ [m]), $m$ the (effective) surge mass/inertia ($17\!\cdot\!10^{6}$ [kg]), and $\{a_{ij}(\cdot),\,b_{j}(\cdot),\,X_u\}$ are (speed-scaled) hydrodynamic derivatives. For these derivatives, dependence on $U$ is mild over moderate-speed cargo-vessel maneuvers; hence we use constant coefficients $a_{ij}(U)\!\equiv\!a_{ij}$, $b_j(U)\!\equiv\!b_j$, with $(\alpha_{11},\alpha_{12},\alpha_{21},\alpha_{22},b_1,b_2,X_u)=(-0.77,-0.34,-3.39,-1.63,0.17,-1.63,83642)$.
%In the baseline case study—large cargo vessels with moderate-speed maneuvers—these derivatives exhibit only weak dependence on $U$; accordingly, we adopt the constant-coefficient approximation $a_{ij}(U)\!\equiv\!a_{ij}$ and $b_{j}(U)\!\equiv\!b_{j}$, with values $\{\alpha_{11},\, \alpha_{12},\, \alpha_{21},\, \alpha_{22},\, b_{1},\, b_{2},\, X_{u},\}=\{-0.77,\,-0.34,\,-3.39,\,-1.63,\,0.17,\,-1.63,\,83642,\}$.

%A central choice in gain scheduling is the scheduler. Two widely used guidelines are that it should (i) vary sufficiently slowly and (ii) capture the dominant plant nonlinearities; under these conditions, scheduled designs tend to inherit local stability/performance across operating conditions \cite{Rugh_20}. Although this is classically illustrated in flight control, the same rationale applies directly to marine autopilots \cite{astrom1994adaptive,aastrom1976identification}.

%For surface vessels, yaw–sway (and yaw–heading) dynamics depend markedly on the through-water surge speed $U$. This variable is directly measurable, evolves on a slow time scale for large ships, and indexes the principal changes of the horizontal-plane maneuvering dynamics. We therefore adopt a MARINER-hull cargo ship as our nonlinear benchmark system \cite{chislett1965planar,fossen1994guidance}.

\subsection{Experimental Setting}\label{sec:expset}
\noindent\textbf{Scheduler and equilibrium-based initialization.} We choose the surge speed $U$ as the scheduling variable ($\rho(k)\!\equiv\!U(k)$): it is directly measurable, varies slowly on large vessels, and strongly shapes yaw–sway dynamics; moreover, $U$ is a state driven by the propeller thrust through the surge model in \eqref{eq:surge}.
%We select the surge speed $U$ as scheduling variable, $\rho(k)\equiv U(k)$. It is directly measurable, evolves on a slow time scale for large ships, and strongly shapes the yaw--sway dynamics; moreover, $U$ appears as a state driven by the propeller thrust $\tau$ via the surge dynamics $\dot U=(-X_u\, U+ \tau)/m$ in \eqref{eq:surge}. The regulated output is the heading $\psi$ with reference signal $\psi_r$.

To motivate the initialization, note that in normal steering the ship executes small deviations about straight-line motion. With the rudder centered, the equilibrium at a given $U_{\mathrm{ss}}$ has $v_{\mathrm{ss}}=r_{\mathrm{ss}}=\delta_{\mathrm{ss}}=0$ and thrust balancing surge resistance, i.e., $\tau_{\text{ss}}=X_u\,U_{\text{ss}}$.
Under the standard slow-variation assumption \cite[Thm.~12.1]{khalil2002nonlinear}, initializing near the equilibrium associated with $\rho(0)=U(0)$ keeps the heading tracking error small and it vanishes as $\rho$ settles. Accordingly, we set $\ve x(0)=[\,0,\,0,\,\psi_r(0)-\kappa,\,U(0)\,]^\top$ and
$\ve u(0)=[\,0,\,X_u U(0)\,]^\top$ with a small bias $\kappa>0$ to avoid the trivial zero-error trajectory. The DeePC windows are seeded by holding $\ve u(0)$ for $T_{\mathrm{ini}}$ samples and using the resulting output measurements for $\ve y_{\mathrm{ini}}$.
Throughout we use $\Delta t=0.1$, $T_{\mathrm{ini}}=10$, and $N=100$.
% To motivate the initialization, note that in normal steering the ship executes small deviations about straight-line motion. With the rudder held at the center position, the steady state consistent with a given surge speed $U_{\text{ss}}$ is $v_{\text{ss}}=0,\; r_{\text{ss}}=0,\; \delta_{\text{ss}}=0,\; \tau_{\text{ss}}=X_u\,U_{\text{ss}}$. Under the discrete-time analogue of Theorem~12.1 in \cite{khalil2002nonlinear}---namely, if the scheduling variable varies slowly and the initial state is sufficiently close to the equilibrium associated with $\rho(0)$, then after a short transient the tracking error $e(k)=\psi(k)-\psi_r(k)$ remains small and tends to zero as the scheduler settles; in particular, if $\rho(k)\to\rho_{\text{ss}}$ and \MAGS{$\dot\rho(k)\to 0$} as $k\to\infty$, then $e(k)\to 0$. In line with this guideline, we initialize near the steady state associated with $ U(0)$, i.e., $
% \ve x(0)=\big[\,0,\,0,\,\psi_r(0)-\kappa,\, U(0)\,\big]^\top,\quad 
% \ve u(0)=\big[\,0,\,X_u\, U(0)\,\big]^\top$, with a small bias $\kappa>0$ to avoid the trivial zero-error trajectory. The DeePC windows are seeded accordingly: $u_{\mathrm{ini}}$ is held at $\ve u(0)$ for $T_{\mathrm{ini}}$ samples, and $\ve y_{\mathrm{ini}}$ is built from the measured outputs produced by this initialization. Throughout we use sampling $\Delta t= 0.1$, $T_{ini}=10$ and $N=100.$\\

\noindent\textbf{Region selection map.}
We realize the selection map $\pi\colon \mathcal{D}_\rho\!\to\!\mathcal{I}$ as a binning of the admissible surge–speed range. Specifically, we define speed bins $[U_i^{\mathrm{lo}},\,U_i^{\mathrm{hi}}]$, $i\in\mathcal{I}$, and set $\pi(U)=i$ whenever $U\in[U_i^{\mathrm{lo}},\,U_i^{\mathrm{hi}}]$. In our experiments we use $M=5$ bins (m/s): $ [0.0,\,2.5],\;\;[2.5,\,5.0],\;\;[5.0,\,9.0],\;\;[9.0,\,12.0],\;\;[11.0,\,14.0]$.
The bin edges reflect hydrodynamic insight and design/sea–trial data (service speed, operating envelope, acceleration capability), guided by classical ship maneuvering/identification sources~\cite{chislett1965planar,aastrom1976identification,astrom1994adaptive}.
%The bin edges are chosen using prior hydrodynamic knowledge and design/sea–trial data (e.g., nominal service speed, operational envelope, acceleration capability), following classical sources on ship maneuvering and identification~\cite{chislett1965planar,aastrom1976identification,astrom1994adaptive}.\\

\noindent\textbf{Data bank and mode slices.}
Because the ship’s states have long time constants, we keep Hankel sizes tractable by partitioning each region’s I/O data into three slices defined by the short-term trend of the regulated output $\psi$: \emph{Up} (rising), \emph{Hold} (approximately constant), and \emph{Down} (falling). At run time, the DeePC predictor for region $i$ uses only the Hankel columns from the slice matching the current trend, preserving locality while reducing the number of columns. In our experiments with $M=5$ regions this yields 15 data subsets; each subset contains 750--1000 samples.\\
\noindent\textbf{DeePC-GS vs.\ baselines computation.}
DeePC-GS uses a bank of $\sim\!15{,}000$ samples—generated by bounded i.i.d.\ uniformly distributed rudder/thrust inputs to ensure PE—\emph{split by region/mode}; online it loads only the \emph{active slice} (typically $n_g\!\approx\!750$–$1000$ Hankel columns). Baselines (no slicing) use one “global” Hankel matrix from a representative $\sim\!6{,}500$ samples spanning all regions/modes. Since the KKT factorization scales as $O(n_g^3)$, this yields $\approx(6500/1000)^3\!=\!275\times$ fewer cubic flops for DeePC-GS (up to $\approx650\times$ if $n_g\!\approx\!750$), with a quadratic memory drop of $\approx42$–$75\times$.

%\emph{DeePC-GS} leverages a bank of $\sim$15{,}000 samples---generated by bounded i.i.d.\ uniformly distributed rudder/thrust inputs to ensure PE---\emph{distributed across region/mode slices} to keep the online computations tractable. At runtime, \emph{DeePC-GS} uses only the predictor corresponding to the \emph{active slice} (typically $n_g\!\approx\!750$--$1000$ Hankel columns), whereas the baselines that do not slice use a single ``global'' Hankel matrix built from a representative subset of $\sim$6{,}500 samples spanning all regions and trend modes. By the $O(n_g^3)$ scaling of the KKT factorization, this yields $(6500/1000)^3\!\approx\!275$ times fewer cubic flops per iteration for \emph{DeePC-GS} (and $\approx\!650$ if $n_g\!\approx\!750$), with a matching quadratic drop in memory ($\sim\!42$--$75$ times smaller).\\
\noindent\textbf{Baseline configurations and input constraints.} We assume full-state measurement, i.e., $\ve y=\ve x=[\,v,\ r,\ \psi,\ U\,]^\top$. 
Inputs are bounded as $\delta\in[-80^\circ,80^\circ]$, $\tau\in[2\times10^3,\,1.5\times10^6]$. 
Unless stated, DeePC-GS uses region-dependent weights $(Q_i,\ve{R}_i,\ve{S}_i,\lambda_{g,i},\lambda_{y,i},\lambda_{u,i})$; baselines (RoKDeePC, Koopman-MPC, Robust DeePC) use the regionwise mean—variations across $i$ are mild so the average is representative. 
Across all runs: $Q=8000$, $(\lambda_g,\lambda_y,\lambda_u)=(100,30,1)$; 
$\ve{R}=\mathrm{diag}(R_\delta,R_\tau)$ with $(R_\delta,R_\tau)=(0.011,\ 4.1\times10^{-4})$; and
$\ve{S}=\mathrm{diag}(\dot R_\delta,\dot R_\tau)$ with $(\dot R_\delta,\dot R_\tau)=(0.48,\ 5.8\times10^{-4})$. 
RoKDeePC is tested with Gaussian, exponential, and degree-5 polynomial kernels; Koopman-MPC uses a third-order polynomial lifting with pairwise products ($n_z=19$).
\subsection{Experimental Results}
We evaluate two ship-steering scenarios with a time-varying heading reference: (i) nominal (no disturbances), benchmarking DeePC-GS against DeePC, RoKDeePC, and Koopman-MPC; and (ii) a following-sea disturbance (stern-to-bow), assessing tracking, effort, and smoothness.

\subsection*{Test $\#$1—nominal tracking}
\begin{figure}[h!]
\centering
    %\vspace{0.4 em}
\includegraphics[trim=0mm 0mm 0mm 0mm,clip, width=1\columnwidth]{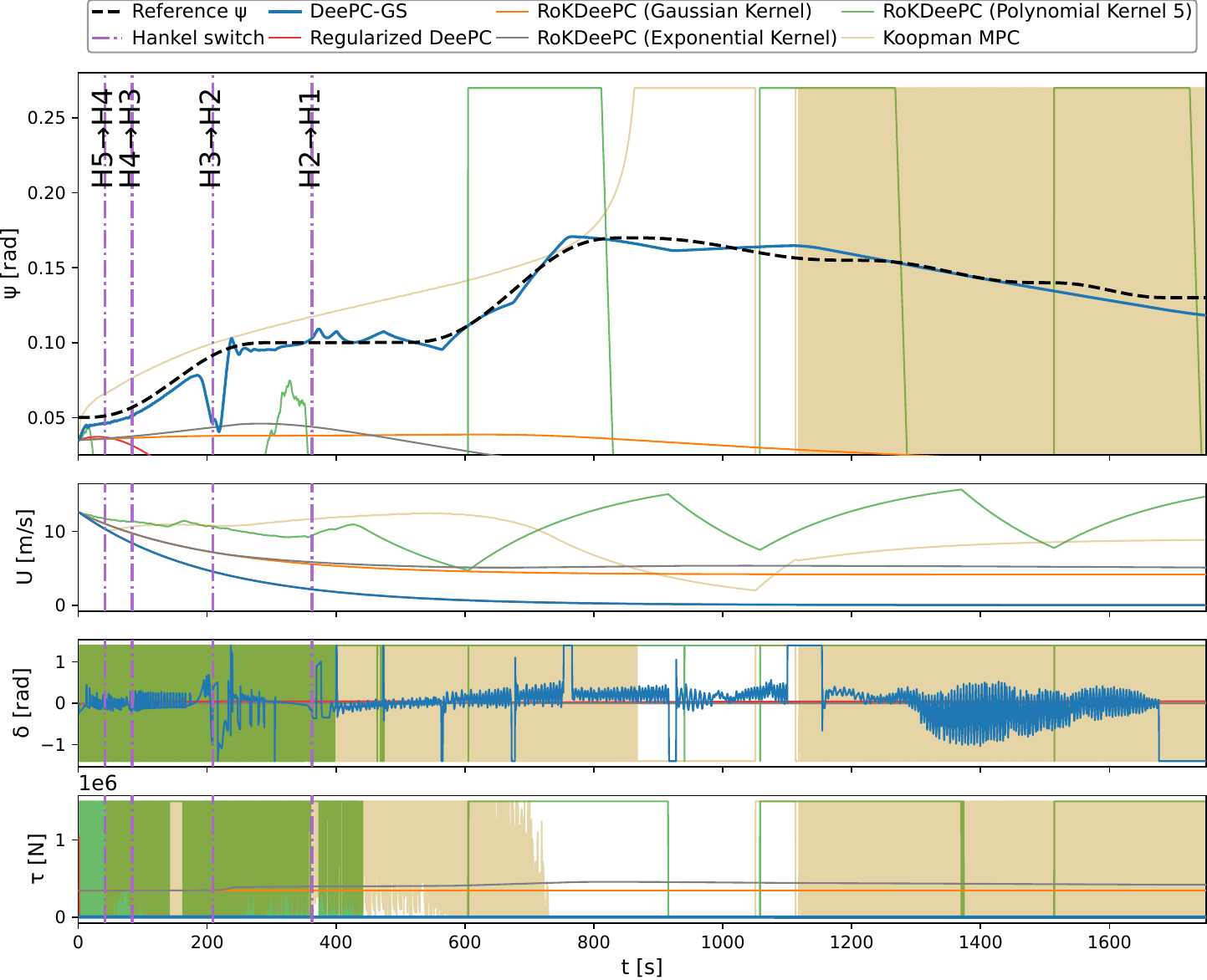}
\caption{Nominal heading tracking. Top: yaw angle $\psi$ (solid) and reference (dashed). Upper middle: surge speed $U$ (scheduling variable). Lower middle: rudder angle $\delta$. Bottom: propeller thrust $\tau$. Vertical purple dashes indicate Hankel switches; $H_j\!\to\!H_k$ labels denote the switch.}
%Vertical purple dashed lines mark Hankel switches; labels $H_j\!\rightarrow\!H_k$ denote the scheduler’s change from local Hankel dataset $H_j$ to $H_k$.}
    \label{fig:tracking}\vspace{-1.15em}
\end{figure}

As shown in Figure~\ref{fig:tracking}, we command the autopilot to follow a time–varying heading reference for $29.2$\,min ($\sim\!30$\,min) to stress long–time robustness on a slow plant. In DeePC-GS the first significant adjustments occur within the initial $\sim$10\,min, during which the surge speed $U$ naturally decays while the optimizer finds feasible solutions with minimal thrust—no terminal constraint is imposed on $U$. A brief tracking dip appears around $200$–$250$\,s as $U$ drifts near the upper boundary of bin $[5,9]~\mathrm{m/s}$, where the local Hankel predictor becomes less representative\footnote{Raising that boundary to, e.g., $6~\mathrm{m/s}$ removes the dip, but we keep it to illustrate that the nonparametric predictor is only locally valid.}; tracking is recovered immediately after switching to the adjacent bin $[2.5,5]~\mathrm{m/s}$. At low speeds the controller increases rudder activity to preserve heading accuracy even when $U\!\approx\!0$. 

Among the baselines, regularized DeePC and the RoKDeePC variants fail to track from early times; the \emph{polynomial} RoKDeePC, in particular, repeatedly drives the actuators into saturation. Koopman-MPC initially tends to follow the reference but, around $t\!\approx\!1100$\,s, develops a high--frequency oscillation, and both $\delta$ and $\tau$ reach their bounds. This behavior is consistent with predictor mismatch: a single global Hankel predictor or a fixed lifted linear model fails to capture the strong surge--speed dependence. Consequently, the controller issues overly aggressive commands that reach actuator limits; once saturated, the closed loop loses phase margin and excites the lightly damped yaw--sway mode, yielding oscillations.
%in particular, repeatedly drives the actuators into saturation. Koopman--MPC initially tends to follow the reference but, around $t\!\approx\!1100$\,s, develops a high--frequency oscillation, and both $\delta$ and $\tau$ hit their bounds. This behavior is consistent with predictor mismatch: a single 
%global Hankel or a fixed lifted lineal model fails to capture the pronounced surge–speed dependence, leading the optimizer to over-command; actuator saturation ensues, the closed loop loses phase margin, and the lightly damped yaw–sway mode is excited. Polynomial kernels further exacerbate the issue via aggressive out-of-manifold extrapolation and ill-conditioned regression, rendering predictions overly sensitive and inducing chatter near constraints.
%global Hankel (regularized DeePC/RoKDeePC) or a fixed lifted linear model (Koopman--MPC) cannot capture the strong surge--speed dependence. Consequently, the controller issues overly aggressive commands that hit the actuator limits; once saturated, the closed loop loses phase margin and excites the lightly damped yaw--sway mode, yielding oscillations. Polynomial kernels exacerbate this: outside the training manifold they extrapolate aggressively, the regression becomes ill--conditioned, and the predicted output turns overly sensitive (chattering near constraints). 
In contrast, DeePC-GS tracks tightly throughout, using regional Hankel slices ($\approx$750–1000 columns) that keep the predictor locally valid, improve conditioning, and reduce computational burden—thus avoiding the above pathologies. Table~\ref{tab:perf_summary} summarizes performance: DeePC-GS leads on all metrics, with solve times reported relative to DeePC-GS (1$\times$).
%In contrast, \emph{DeePC-GS} maintains tight tracking throughout. At each step it uses regional Hankel slices ($\approx$750--1000 columns), preserving predictor validity, improving numerical conditioning, and reducing computational burden—thereby avoiding the aforementioned pathologies. A quantitative summary in Table~\eqref{tab:perf_summary} shows \emph{DeePC–GS} outperforming the alternatives across all metrics; solve times are shown relative to \emph{DeePC-GS} (1$\times$).
\begin{table}[h!]
%\vspace{-1em}
    \centering
    \caption{RMSE, input energy, and normalized solve time.}
    \label{tab:perf_summary}
    \scriptsize
    \setlength{\tabcolsep}{4pt}
    \begin{tabular}{l|c|c|c}
        \hline
        Method & RMSE$_\psi$ & Input Energy $E_{\mathrm{tot}}$ & Time \\
        \hline
        DeePC-GS & $0.0075619$ & $8.806\times10^{9}$  & $1\times$ \\
        Regularized DeePC & $1.4445$ & $1.126\times10^{12}$ & $6\times$ \\
        RoKDeePC (Gaussian) & $0.10361$ & $2.114\times10^{14}$ & $16\times$ \\
        RoKDeePC (Exponential) & $0.11090$ & $3.086\times10^{14}$ & $16\times$ \\
        RoKDeePC (Polynomial-5) & $1.1754$ & $2.464\times10^{15}$ & $16\times$ \\
        Koopman MPC & $1.5171$ & $1.739\times10^{15}$ & $1.1\times$ \\
        \hline
    \end{tabular}
    \vspace{-2em}
\end{table}
%Among baselines, regularized DeePC and RoKDeePC approaches fail to track from early times, whereas Koopman–MPC achieves reasonable tracking albeit with higher control effort/rate. At each step, \emph{DeePC–GS} uses regional Hankel slices ($\approx$750–1000 samples), whereas the \emph{DeePC–type} baselines  rely on a single global Hankel ($\approx$6500 samples); Koopman–MPC instead uses a fixed lifted linear model. Accordingly, local data selection improves conditioning and reduces computation without loss of performance; hence robust–DeePC’s failure is expected—one global Hankel cannot capture parameter–varying ship dynamics, so an LTI predictor lacks local validity across regimes. A quantitative summary (RMSE in $\psi$, input energy, and input–rate cost) in Table~2 shows \emph{DeePC–GS} outperforming the alternatives across all metrics.
\subsection*{ Test $\#$2—following-sea disturbance}
\begin{figure}[h!]
\centering
    \vspace{0.5 em}
    %\vspace{0.4 em}
\includegraphics[trim=0mm 0mm 0mm 0mm,clip, width=1\columnwidth]{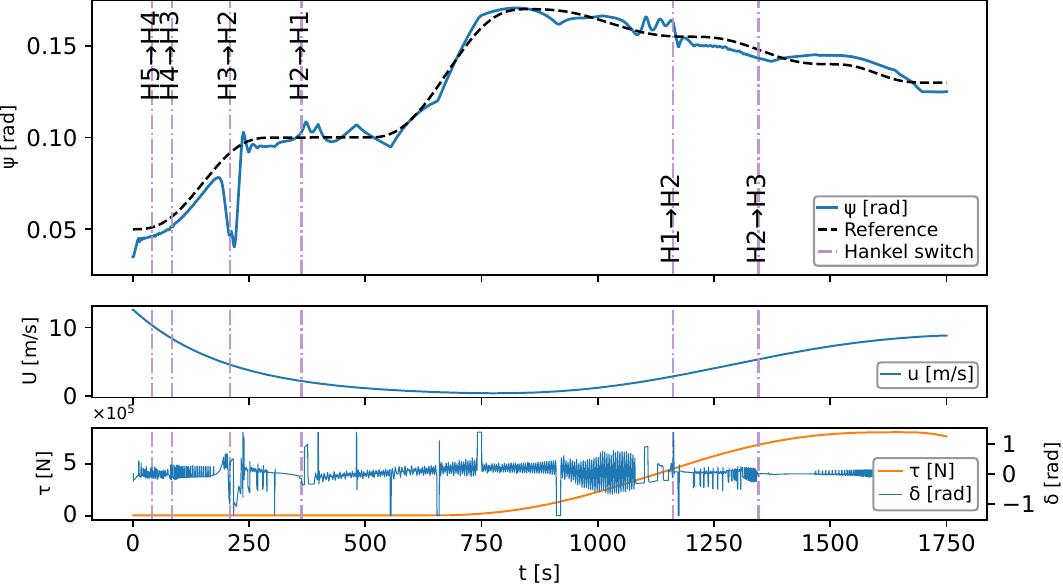}
\caption{Heading tracking with disturbance. Top: yaw angle $\psi$ (solid) and reference (dashed). Middle: surge speed $U$ (scheduling variable). Bottom: propeller thrust $\tau$ (left axis) and rudder angle $\delta$ (right axis). Vertical purple dashes indicate Hankel switches; $H_j\!\to\!H_k$ labels denote the switch.}
    \label{fig:disturbance} \vspace{-1.4em}
\end{figure}

A stern--to--bow wave is applied, ramping in at $t\approx 600$\,s, implemented as an additive actuator disturbance so the applied thrust satisfies $\tau_{\mathrm{applied}}=\tau_{\mathrm{control}}+\tau_{\mathrm{wave}}$. As shown in Figure~\ref{fig:disturbance}, the vessel then returns to the $[5,9]$\,m/s bin with a smooth switch while maintaining tight heading tracking. At speeds \emph{above} bin~1 ($[0,2.5]$\,m/s), tracking improves, owing to the reduced control authority required to drive $\psi$ to the reference compared to very low speeds. This aligns with the MARINER hull specification: the nominal speed ($\approx 7.7$\,m/s) lies near the center of the operating range.
%A stern--to--bow wave is applied, ramping in at $t\approx 600$\,s; the vessel then returns to the $[5,9]$\,m/s bin with a smooth switch while maintaining tight heading tracking. At speeds \emph{above} bin~1 ($[0,2.5]$\,m/s), tracking improves owing to the reduced control authority required to drive $\psi$ to the reference compared with very low speeds. This aligns with the MARINE hull specification: the nominal speed ($\approx 7.7$\,m/s) lies near the center of the operating range.
%focusing on the accuracy of the estimated structured singular value and the impact of the underlying uncertainty structure $\Delta$ on its performance.

\section{Conclusion} \label{sec: conclusion}
We presented \emph{DeePC-GS}, a gain-scheduled, fully data-driven predictive controller that switches among local Hankel predictors via a measurable scheduler. On a nonlinear ship-steering benchmark, \emph{DeePC-GS} achieves accurate tracking with tractable computation and outperforms state-of-the-art data-driven MPC baselines. Future work includes learning operating regions and the scheduling map directly from data, removing \emph{a priori} knowledge about the partitioning.
%In this paper, we have introduced \emph{DeePC–GS}, a gain-scheduled, fully data-driven predictive control scheme that leverages local Hankel predictors selected online via a measurable scheduling variable. Numerical simulations on a nonlinear ship-steering benchmark demonstrate that \emph{DeePC–GS} achieves accurate tracking with tractable computation and outperforms state-of-the-art data-driven MPC baselines. As future work, we aim to remove \emph{a priori} region specifications by automatically discovering operating regions and the scheduling map directly from data.
%Notice that this transformation does not change the objective of estimating the unknown parameter $\theta = (a \ b)^T$. So, in order to estimate $\theta$ corresponding to the stochastic volatility model~\eqref{eq: sim_dynamic_2} using TSGBM, we consider~\eqref{eq: sim_dynamic_2_eq} instead. 

%To illustrate the performance of minimax TSGBM, we first show 

\bibliography{References}

%%%%%%%%%%%%%%%%%%%%%%%%%%%%%%%%%%%%%%%%%%%%%%%%%%%%%%%%%%%%%%%%%%%%%%%%%%%%%%%%

%%%%%%%%%%%%%%%%%%%%%%%%%%%%%%%%%%%%%%%%%%%%%%%%%%%%%%%%%%%%%%%%%%%%%%%%%%%%%%%%

\end{document}